
\input amstex

\documentstyle{amsppt}

\loadbold

\magnification=\magstep1

\pageheight{9.0truein}
\pagewidth{6.5truein}

\NoBlackBoxes

\def\goi{\frak{g}_0 \oplus \frak{g}_1}
\def\g{\frak{g}}
\def\go{\frak{g}_0}
\def\gi{\frak{g}_1}

\def\c{\frak{c}}
\def\co{\frak{c}_0}

\def\GKdim{\operatorname{GKdim}}

\def\a{\frak{a}}
\def\ao{\frak{a}_0}
\def\ai{\frak{a}_1}

\def\b{\frak{b}}

\def\U{{\overline U}}
\def\ad{\operatorname{ad}}
\def\Z{\Bbb{Z}}
\def\ad{\operatorname{ad}}
\def\pl{\operatorname{pl}}

\def\SuperVect{\text{SuperVect}}

\def\AndEtiGel{1}
\def\Beh{2}
\def\Bro{3}
\def\BroGoo{4}
\def\Fis{5}
\def\GooLet{6}
\def\Kos{7}
\def\KraLen{8}
\def\McCRob{9}
\def\Maj{10}
\def\Mon{11}
\def\Rad{12}
\def\Sch{13}

\topmatter

\title An affine PI Hopf algebra not finite over a normal
commutative Hopf subalgebra \endtitle

\date December 2001. \enddate

\rightheadtext{An affine PI Hopf algebra}

\author Shlomo Gelaki and Edward S. Letzter \endauthor

\abstract In formulating a generalized framework to study certain
noncommutative algebras naturally arising in representation theory,
K. A. Brown asked if every finitely generated Hopf algebra satisfying a
polynomial identity was finite over a normal commutative Hopf subalgebra. In
this note we show that Radford's biproduct, applied to the enveloping algebra
of the Lie superalgebra $\pl(1,1)$, provides a noetherian prime
counterexample. \endabstract

\address Department of Mathematics, Technion-Israel Institute of
Technology, Haifa 32000, \linebreak Israel \endaddress

\email gelaki\@math.technion.ac.il \endemail

\address Department of Mathematics, Temple University, Philadelphia,
PA 19122, USA \endaddress

\email letzter\@math.temple.edu \endemail

\thanks The first author's research was supported by the Technion
V.P.R. Fund - Loewengart research fund and by the Fund for the
Promotion of Research at the Technion. The second author's
research was supported in part by NSF grant DMS-9970413. This
research was begun during the second author's visit to the
Technion in August 2001, and he is grateful to the Technion for
its hospitality. \endthanks

\endtopmatter

\document

\head 1. Introduction \endhead

Noncommutative affine PI Hopf algebras arise in fundamental
representation-theoretic contexts: as enveloping algebras of
finite dimensional restricted Lie algebras in positive
characteristic, as group algebras of finitely generated
abelian-by-finite groups, and as quantizations of algebraic groups
-- and their corresponding Lie algebras -- at roots of unity. (By
{\sl affine\/} we mean ``finitely generated over the ground field
as an associative algebra,'' and by PI we mean ``satisfying a
polynomial identity.'') In all of the preceding examples, the Hopf
algebra of interest is finitely generated as a module over a
normal affine commutative Hopf subalgebra. Therefore, in
developing a general structure theory to handle these examples, K.
A. Brown asked \cite{\Bro, Question B} whether every affine PI
Hopf algebra was finitely generated as a left (or right) module
over some normal commutative Hopf subalgebra (cf\.
\cite{\BroGoo}). In this note we present a counterexample; we
describe a prime affine noetherian PI Hopf algebra not finitely
generated as a left or right module over any of its normal
commutative Hopf subalgebras. This counterexample also arises in a
representation-theoretic context, resulting from an application of
Radford's biproduct \cite{\Rad} to the enveloping algebra of the
Lie superalgebra $\pl(1,1)$. (In the terminology of \cite{\Maj},
the counterexample is the ``bosonization'' of the enveloping
algebra of $\pl(1,1)$.)

\head 2. Preparatory Results \endhead

Let $k$ be an algebraically closed field of characteristic zero. We
assume throughout this note that $k$ is the ground field for all of
the vector spaces mentioned (explicitly or implicitly). In particular,
all algebras (associative and otherwise) are over $k$. All of the
associative algebras discussed below will be assumed to be unital.

The reader is referred to \cite{\McCRob} for background on
noncommutative rings, to \cite{\Mon} for background on Hopf algebras,
and to \cite{\Sch} for background on Lie superalgebras.

\subhead 2.1 \endsubhead (i) For all of the Hopf algebras $H$ we will
encounter, we will use $\Delta_H$ to denote the coproduct,
$\epsilon_H$ to denote the counit, and $S_H$ to denote the
antipode. The subscripts will be dropped when the meaning is
clear.

(ii) Let $H$ be a Hopf algebra. Recall the left and right adjoint
actions of $H$ on itself,
$$(\ad _l a)(b) = \sum a_1 b S(a_2), \quad \text{and} \quad (\ad _r
a)(b) = \sum S(a_1)b a_2 ,$$
for $a, b \in H$, where $\Delta (a) = \sum a_1 \otimes a_2$. Also
recall that a Hopf subalgebra $A$ of $H$ is {\sl normal\/} in $H$
provided $(\ad_lH)(A)$ and $(\ad _r H)(A)$ are contained in $A$.

\subhead 2.2 \endsubhead We next briefly review some superalgebra.

(i) A vector space $V$ equipped with a $\Z_2$-grading is termed a {\sl
supervector space\/}; we indicate the grading by writing $V = V_0
\oplus V_1$. The vectors in $V_0$ are called {\sl even\/} and those in
$V_1$ are called {\sl odd}. We define the parity of a homogeneous
vector $v \in V$ to be $p(v)=0$ if $v$ is even and $p(v)=1$ if $v$ is
odd.

We will use $\SuperVect$ to denote the category of supervector spaces
with homogeneous morphisms. For $V, W \in \SuperVect$ the
commutativity constraint $V\otimes W\to W\otimes V$ is different from
the classical one and is given by the formula
$$v\otimes w\mapsto (-1)^{p(v)p(w)}w\otimes v,$$
for homogeneous $v \in V$ and $w \in W$.

(ii) One can define Lie algebras, bialgebras, Hopf algebras, etc\. in
$\SuperVect$ by replacing the ordinary structure maps with morphisms
in $\SuperVect$. The resulting objects are referred to, respectively,
as Lie superalgebras, superbialgebras, Hopf superalgebras, etc.

For example, a Hopf superalgebra $H$ is a $\Z_2$-graded associative
algebra, possessing a coassociative multiplicative (in the
super-sense) morphism $\Delta\colon H\rightarrow H\otimes H$ in
$\text{SuperVect}$, and further equipped with a counit $\epsilon$ and
antipode $S$ satisfying the expected axioms. By ``multiplicativity in
the super-sense'' we mean that $\Delta$ satisfies the relation
$$\Delta(ab)=\sum (-1)^{p(a_2)p(b_1)}a_1b_1\otimes a_2b_2$$
for all homogeneous $a,b\in H$ (where $\Delta(a)=\sum a_1\otimes
a_2$, $\Delta(b)=\sum b_1\otimes b_2$ and the components are
homogeneous). This formulation is required because the tensor
product of algebras $A$ and $B$ in $\text{SuperVect}$ is defined
to be the vector space $A\otimes B$ equipped with the
multiplication
$$(a\otimes b)(a'\otimes b'):=(-1)^{p(a')p(b)}a a'\otimes bb',$$
for homogeneous $a, a' \in A$ and $b, b' \in B$. (Hopf
superalgebras are referred to as ``graded Hopf algebras'' in
\cite{\Kos}.)

(iii) Henceforth, ``graded'' will mean ``$\Z_2$-graded.''

\subhead 2.3 \endsubhead (i) The Hopf superalgebras we are most concerned
with in this note arise as follows: Let $\g$ be a Lie superalgebra,
and let $U(\g)$ be the enveloping algebra of $\g$. Then $U(\g)$ is a
Hopf superalgebra with coproduct $\Delta$, counit $\epsilon$, and
antipode $S$ determined by
$$\Delta (x) = x\otimes 1 + 1 \otimes x, \quad S(x) = -x, \quad
\epsilon(x) = 0,$$
for all $x \in \g$. Note that $U(\g)$ is supercocommutative with
no non-trivial grouplike elements.

(ii) Now suppose that $U'$ is a sub Hopf superalgebra of $U(\g)$.  By
a theorem of Kostant \cite{\Kos, Theorem 3} (cf\.  \cite{\AndEtiGel}),
there exists a sub Lie superalgebra $\a$ of $\g$ such that $U' = U(\a)
\subseteq U$. (We can apply Kostant's theorem because $U$, and hence
also $U'$, is supercocommutative with no non-trivial grouplike
elements). Note that $\ao$ is contained in $\go$ and that $\ai$ is
contained in $\gi$.

\subhead 2.4 \endsubhead Let $U$ be any Hopf superalgebra.

(i) Recall (see, e.g., \cite{\AndEtiGel}) that Radford's biproduct
\cite{\Rad} allows us to associate to $U$ a Hopf algebra $\U$ in the
following way. As an algebra,
$$\U=U\# k[\langle t \rangle] \cong U\# k[\Z/2\Z ]$$
is the smash product of $U$ and $k[\langle t \rangle]$, where $t$ is
the graded algebra automorphism of $U$ (of order $2$) acting on $U$ by
parity (i.e. $tat = (-1)^{p(a)}a$ for all homogeneous $a \in U$). The
coproduct $\Delta$, counit $\epsilon$, and antipode $S$ of $\U$ are
determined by
$$\Delta (x)=\sum x_1\#t^{p(x_2)}\otimes x_2\# 1, \quad \Delta(t) =
t\otimes t,$$
and
$$S(x) =t^{p(x)} S_U(x),\quad \epsilon(x)=\epsilon_U(x),\quad S(t) =
t,\quad \epsilon(t) = 1 ,$$
where $x\in U$ is homogeneous and $\Delta_U(x) = \sum x_1 \otimes
x_2$ with homogeneous components.

(ii) In particular, $K:=k[\langle t \rangle ]$ is a Hopf subalgebra of
$\U$, and there is a unique Hopf algebra projection
$$\U @> \quad \pi \quad >> K ,$$
sending $u\#h$ to $\epsilon(u)h.$ Observe that
$$U = \{ x \in \U : (\text{id}\otimes \pi)\Delta(x) = x\otimes 1 \},$$
where $\text{id}$ denotes the identity map, is the algebra of
$K$-coinvariants of $\U$ (see \cite {\Rad}).

(iii) Also note that $U$ is a subalgebra of $\U$, and there is a
unique coalgebra projection
$$\U @> \quad \Pi \quad >> U,$$
sending $u\#h$ to $u\epsilon(h).$

(iv) Lastly, by declaring $t$ to be even, we can extend the grading on $U$
to a grading on $\U$. Note that $ta = (-1)^{p(a)}at$ for all
homogeneous $a \in \U$.

\subhead 2.5 \endsubhead Let $\g$ be a Lie superalgebra with
enveloping algebra $U := U(\g)$.

(i) The Hopf algebra structure on $\U$ follows from
$$\Delta (x) = x\otimes 1 + t^{p(x)} \otimes x, \; S(x) =-t^{p(x)} x,
\; \epsilon(x) = 0, \; \Delta(t) = t\otimes t, \; S(t) = t, \;
\epsilon(t) = 1 ,$$
for homogeneous $x \in \g$ (see, e.g., \cite{\Fis}).

(ii) Now choose homogeneous $a,b \in \g$. Write $\Delta_U(a) = \sum
a_1 \otimes a_2$. Then
$$\multline (\ad _\g a)(b) = [a,b] = ab - (-1)^{p(b)p(a)}ba = \sum a_1
b (-1)^{p(b)p(a_2)}S_U(a_2) = \\ \sum (a_1\#1)(b\#1)
(-1)^{p(b)p(a_2)}(1\#t^{p(a_2)})S(a_2\#1) = \\ \sum (a_1\#
t^{p(a_2)})(b\#1)S(a_2\#1) = (\ad_l a)(b) , \endmultline$$
where $\ad_l$ denotes the left adjoint action of $\U$ on itself,
and $\ad_\g$ denotes the adjoint (graded) representation of $\g$
in itself. Consequently, if $A$ is a normal Hopf subalgebra of
$\U$ such that $A \cap U = U(\a)$, for some Lie superalgebra $\a$
of $\g$, then $\a$ is an ideal of $\g$.

(iii) Note that $U(\go)$ is both a sub Hopf superalgebra of $U$ and a Hopf
subalgebra of $\U$.

\proclaim{2.6 Lemma} Let $U$ be any Hopf superalgebra,
and let $A$ be a Hopf subalgebra of $\U$ containing $t$. Then

{\rm (i)} $A \cap U$ is a sub Hopf superalgebra of $U$. In particular,
if $U$ is the enveloping algebra of a Lie superalgebra $\g$, then $A
\cap U = k\langle a \rangle = U(\a)$ for some sub Lie superalgebra
$\a$ of $\g$.

{\rm (ii)} $A = (A\cap U)\# K$.

\endproclaim

\demo{Proof} Set $\Delta := \Delta_\U$ and $S := S_\U$. To start,
$\pi(A) = K$ since $t \in A$.  Also, the composition of Hopf algebra
homomorphisms
$$K @> i >> A @> i >> \U @> \pi >> K ,$$
where $i$ denotes the inclusion map, reduces to the identity map on
$K$. Hence, by Theorem 3 in \cite{\Rad}, we obtain a biproduct
decomposition
$$A = \tilde A \# K,$$
where, by construction, $\tilde A=\Pi(A)$ is a Hopf superalgebra.
Therefore, it is straightforward to verify that
$$\tilde A = \{ x \in A : (\text{id}\otimes \pi)\Delta(x) = x\otimes 1
\} = A \cap U,$$
is a subalgebra of $U$. Moreover, since $A$ and $U$ are stable under
the (left or, equivalently, right) adjoint action of $t$, we see that
$A \cap U$ is stable under the (left or, equivalently right) adjoint
action of $t$. Therefore, $A\cap U$ is a graded subalgebra of $U$.
Also, since $\Pi$ is a coalgebra map, we can conclude that $A \cap U =
\tilde A = \Pi(A)$ is a subcoalgebra of $U$.

Finally, if $a$ is a homogeneous element of $A\cap
U$, then $S(a)= t^{p(a)}S_U(a) \in A$. Hence $S_U(a) = t^{p(a)}S(a) \in A$,
and we see that $S_U(A\cap U) \subseteq A\cap U$. In view of (2.3ii),
the lemma now follows. \qed\enddemo

\head 3. The Example \endhead

\subhead 3.1 \endsubhead In this section we let $\g=\goi$ denote the
Lie superalgebra $\pl(1,1)$ of $2{\times}2$ matrices. Fix a basis for
$\g$,
$$x = \bmatrix 1 & 0\\ 0 & 1 \endbmatrix, \quad y = \bmatrix 1 & 0 \\
0 & 0 \endbmatrix, \quad u = \bmatrix 0 & 1 \\ 0 & 0 \endbmatrix,
\quad v = \bmatrix 0 & 0\\ 1 & 0\endbmatrix .$$
Recall that $\go$ is generated by $x$ and $y$, while $\gi$ is
generated by $u$ and $v$. Set $U := U(\g)$. As in (2.4i), let $\U
:= U \# k[\langle t \rangle ]$. The relations defining $\U$ as an
algebra are
$$\multline xh = hx \quad \text{for all $h \in \U$}, \quad yu - uy = u, \quad
yv -vy = -v, \quad uv + vu = x, \\ \quad u^2 = v^2 = 0, \quad tx - xt = ty -
yt = 0, \quad tu = -ut, \quad tv = -vt, \quad t^2 = 1 . \endmultline $$

We can now state the main result of this paper.

\proclaim{3.2 Theorem} The Hopf algebra $\U$ is a
noetherian prime affine PI Hopf algebra not finitely generated, as a
right or left module, over any of its normal commutative Hopf
subalgebras.  \endproclaim

The remainder of this section is devoted to a proof of the theorem

\subhead 3.3 \endsubhead (i) Note that $u$ is an $\ad_\g(y)$-eigenvector
with eigenvalue $1$ and that $v$ is an $\ad_\g(y)$-eigenvector with
eigenvalue $-1$.

(ii) Let $w$ be an $\ad_\g(y)$-eigenvector of $\gi$. By (i), the
corresponding eigenvalue is $\pm 1$. Hence, in $\U$, for all
positive integers $n$, either $y^nw=w(y+1)^n$ or $y^nw=w(y-1)^n$.

\subhead 3.4 \endsubhead Let $F$ be an algebra containing a subalgebra
$E$, and suppose that $F$ is finitely generated on the right (or left)
as an $E$-module. Noting that $F$ is a ring subquotient of a matrix
ring over $E$, we see that $E$ is PI if and only if $F$ is PI. Next,
we will also require some basic and well-known properties of
Gelfand-Kirillov dimension (abbreviated to ``GK dimension'' or
``GKdim''); see \cite{\KraLen} or \cite{\McCRob, Chapter 8} for
background. In particular, $\GKdim E = \GKdim F$ (see, e.g.,
\cite{\KraLen, 5.5}). The GK dimension of the commutative polynomial
algebra $k[x_1,\ldots,x_n]$ is equal to $n$, and the enveloping
algebra of an $n$-dimensional Lie algebra has GK dimension equal to
$n$.

\subhead 3.5 \endsubhead Note that $k[x,y] = k\langle x, y \rangle
\subset \U$ is a commutative polynomial ring in $x$ and $y$.
Moreover, $\U$ is finitely generated as a left and a right $k[x,
y]$-module. Therefore, by (3.4), $\U$ is an affine noetherian PI
algebra of GK dimension $2$.

\subhead 3.6 \endsubhead Suppose that $\a$ is a sub Lie superalgebra of $\g$
such that $\ao \ne \go$. Then $\dim \ao \leq 1$, and it follows that $\GKdim
U(\ao) \leq 1$. It follows from the Poincar\'e-Birkhoff-Witt Theorem for Lie
superalgebras (PBW Theorem) that $U(\a)$ is finitely generated as a left
$U(\ao)$-module. Hence, by (3.4), we see that $\GKdim U(\a) \leq 1$.

\subhead 3.7 \endsubhead Let $C$ be a commutative normal Hopf
subalgebra of $\U$. We now prove that $\U$ cannot be finitely
generated as a left or right module over $C$.

\subhead Case 1 \endsubhead Assume that $t \not\in C$.

Set $A := C\langle t \rangle$, and observe that $A$ is a Hopf subalgebra of
$\U$. Since $C$ is stable under the left and right adjoint actions of $t$, we
see that $C$ is a graded associative subalgebra of $\U$. Since $t$ is
homogeneous, the grading on $C$ extends to a grading on $A$.  For $a,b \in W
:= C_0 \cup C_1 \cup \{ t \}$, we further see that $ab = \pm ba$. Also, $A$ is
generated as a $k$-algebra by $W$.

Now choose $a_1,\ldots,a_n \in W$. Set $A' := k\langle a_1,\ldots,a_n
\rangle$, and set $Z' := k\langle a_1^2,\ldots, a_n^2 \rangle$. Then
$Z'$ is central in $A$, and $A'$ is generated as a left (or right)
$Z'$-module by the products $a^{i_1}_1 \cdots a^{i_n}_n$, for
$i_1,\ldots,i_n \in \{0,1\}$. In particular, $A'$ is finitely
generated as a left (or right) module over the commutative noetherian
ring $Z'$. It follows that $A$ is integral over its center, in the
sense of \cite{\McCRob, 5.3.2}.

Since $t \in A$, it follows from (2.6i) that $A\cap U = U(\a)$ for some sub
Lie superalgebra $\a$ of $\g$.

\subsubhead Subcase 1a \endsubsubhead Assume that $\ao = \go$ (and so
$y \in A$).

First suppose that $\ai \ne 0$, ensuring that $\ai$ contains an
$\ad_\g(y)$-eigenvector. We can impose a non-negative $\Z$-grading
on $\g$ by setting $\deg(y)=0$, $\deg(u)=\deg(v)=1$, and $\deg(x)
= 2$. This $\Z$-grading restricts to $\a$ since $\a$ must be
generated by $\Z$-homogeneous elements. Moreover, by declaring
$\deg(t) = 0$, we obtain a $\Z$-grading on $A$. Note that the
$0$-degree $\Z$-homogeneous component of $A$ is $k[y,t]$, and also
note that the $\Z$-grading on $A$ restricts to the center of $A$.
It now follows from the PBW Theorem, and (3.3ii), that the
$0$-degree $\Z$-homogeneous component of the center is $k$.
However, we learned above that $y$ must be integral over the
center of $A$, and so $y$ must be integral over the $0$-degree
$\Z$-homogeneous component of the center of $A$. Because $y$
cannot be integral over $k$, we have arrived at a contradiction.
Therefore, $\ai = 0$.

It now follows that $A \cap U = U(\go)$, and it follows from (2.6ii)
that $A = U(\go)\#K$. In particular, $A$ is a pointed cocommutative
Hopf algebra, and the set of grouplike elements of $A$ is $\{1,
t\}$. Also, $\go$ is the Lie algebra of primitive elements of
$A$. Now, since $t \not\in C$, we see that $C$ is a cocommutative
connected Hopf algebra, and again by a theorem of Kostant (see, e.g.,
\cite{\Mon, 5.6.5}), we find that $C = U(\c)$, where $\c$ is the Lie
algebra of primitive elements of $C$. Of course, $\c \subseteq \go$,
and so $\c = \co$ is an even sub Lie superalgebra of $\g$. Moreover,
it follows from (2.5ii) that $\c$ is an ideal of $\g$. Now, any ideal
of $\g$ containing $y$ must also contain $u$ and $v$, and so we see
that $y \not\in \c$. Thus $\c \ne \go$. Therefore,
$\GKdim C \leq 1$, by (3.6), and so $\U$ cannot be a finitely
generated left or right $C$-module, since $\GKdim\U = 2$.

\subsubhead Subcase 1b \endsubsubhead Assume that $\ao \ne \go$.

It now follows from (3.6) that $\GKdim A \cap U \leq 1$. By
(2.6ii), $A$ is finitely generated as a left $A\cap U$-module, and
so by (3.4), $\GKdim A \leq 1$. Therefore, $\U$ cannot be finitely
generated as a left or right $A$-module. Hence $\U$ cannot be
finitely generated as a left or right $C$-module.

\subhead Case 2 \endsubhead Assume that $t \in C$.

By (2.6i), it follows that $C \cap U = U(\c)$, for some sub Lie
superalgebra $\c$ of $U$. Moreover, it follows from (2.5ii) that
$\c$ is an ideal of $\g$. Consequently, if $\co = \go$, then $\c =
\g$. But $\c \ne \g$, since otherwise $U$ is equal to the
commutative algebra $C\cap U$.  So $\co \ne \go$, and we deduce
from (3.6) that $\GKdim C\cap U \leq 1$.  Hence, by (2.6ii) and
(3.4), $\GKdim C \leq 1$. We conclude that $\U$ cannot be a
finitely generated left or right $C$-module.

\subhead 3.8 \endsubhead It only remains to prove that $\U$ is
prime. Let $R$ be a prime noetherian ring, and let $\tau$ be an
automorphism of $R$ of order $2$. Let $T$ denote the Ore extension
$R[Y; \tau]$. It follows, for example, from \cite{\GooLet, 2.3i--iii}
that $T/\langle Y^2 - 1 \rangle$ is prime. As noted in \cite{\Beh},
$U$ is prime. Therefore, since $\U$ is isomorphic as a ring to $U[Y;
t]/\langle Y^2 - 1 \rangle$, we see that $\U$ is prime.

\remark{{\bf 3.9} Remark} Let $\b$ be the sub Lie superalgebra of $\g$
generated by $y$ and $u$. Using the methods of this note, one can show
that $U(\b)\# k[\langle t \rangle]$ is a noetherian Hopf PI algebra,
of GK dimension $1$, not finitely generated as a right or left module
over any of its normal commutative Hopf subalgebras. However, the
ideal $\langle u \rangle$ of $U(\b)\# k[\langle t \rangle]$ is
nilpotent, and so $U(\b)\# k[\langle t \rangle]$ is not semiprime. We
therefore ask: Does there exist a prime affine PI Hopf algebra of GK
dimension $1$ that is not finitely generated as a left or right module
over any of its normal commutative Hopf subalgebras?  \endremark

\head 4. Acknowledgment \endhead

We are grateful to K. A. Brown and Pavel Etingof for reading
preliminary versions of this manuscript, and for numerous useful and
stimulating conversations.

\newpage

\enddocument